\documentclass{article}

\usepackage{enumerate}
\usepackage{amsmath,amssymb}
\usepackage{ascmac}
\usepackage{amscd}

\usepackage[top=30truemm,bottom=30truemm,left=40truemm,right=40truemm]{geometry}

\newtheorem{theorem}{Theorem}[section]
\newtheorem{corollary}{Corollary}[section]
\newtheorem{Prop}{Proposition}[section]
\newtheorem{definition}{Definition}[section]

\newtheorem{lemma}{Lemma}[section]
\newtheorem{Exa}{Example}[section]
\newtheorem{remark}{Remark}[section]
\newtheorem{proof}{Proof.}

\numberwithin{equation}{section}

\title{Digital pattern and transcendence via\\ generalized $k$-regular sequences }

\author{Eiji Miyanohara \\
E-mail: j1o9t5acrmo@fuji.waseda.jp}
\date{\today}

\begin{document}

\date{}
\maketitle
\begin{abstract}
In this paper, we prove that there are uncountable many real transcendental numbers, which are  generated by digital pattern sequences. This generalizes the main theorem in Morton and Mourant, which states the existence of countable many similar numbers. Our method relies on the combinatorial quantitative transcendence criterion established by Adamczewski and Bugeaud and properties of generalized $k$-regular sequences, which is introduced by this paper.
\newline
\newline

\end{abstract}

\section{Introduction}\label{s:1}
\quad
First, we introduce the Rudin-Shapiro sequence by digital pattern counting.
Let $k $ be an integer with $2\le k$. We define the base $k$-representation of the non-negative integer $n$ as follows
\begin{align}\label{eq1.1}
n=\sum_{q=0}  s_{n,q} k^{q},
\end{align}
where $0\le s_{n,q}\le k-1$. Let $A$ be a set and ${A}^{*}$ be the free monoid generated by $A$. For a positive integer $d$, we denote by $P:=p_1 \cdots p_d \in {\{0,1,\cdots k-1\}}^{*}/{\{0\}}^{*}$ a pattern. For a pattern $P$, $e_P(n)$ denote the number of occurrences of $P$ in the base $k$-representation of $n$. (See the definition of pattern sequences \cite{AlS2},\cite{MM},\cite{ST})For an integer $L$ with $2\le L$, we define the sequence
${(e^L_P(n))}_{n\ge0}$ by
\begin{align}
e^L_P(n)\equiv e_P(n) \pmod L
\end{align}
where $0\le e^L_P(n)\le L-1$. Then ${(e^2_{11}(n))}_{n\ge0}$, where $k=L=2$, is known as the Rudin-Shapiro sequence.
The Rudin-Shapiro sequence is a famous automatic sequence.
Morton and Mourant \cite{MM}
 and Adamczewski and Bugeaud \cite{AB1} have proved the following result.
\begin{theorem}\label{MM-AB}{\rm(\cite{MM},\cite{AB1})}
{\rm
Let $\beta \ge L$ be an integer. Then $\sum_{n=0}^\infty \frac{e^L_P(n)}{ {\beta}^{n+1}}$ is a transcendental number unless $k=3, P=1$ and $L=2$.}
\end{theorem}
The proof of Theorem \ref{MM-AB} relies on the periodicity of ${(e^L_P(n))}_{n\ge0}$\cite{MM} and the Cobham conjecture that was settled by Adamczedwski and Bugeaud (Theorem $2$ in \cite{AB2}). More precisely, Morton and Mourant proved that ${(e^L_P(n))}_{n\ge0}$ is a $k$-automatic sequence.
Furthermore, they proved that ${(e^L_P(n))}_{n\ge0}$ is periodic if and only if ${(e^L_P(n))}_{n\ge0}$ is purely periodic, which enabled them to prove that ${(e^L_P(n))}_{n\ge0}$ is periodic if and only if ${(e^L_P(n))}_{n\ge0}$ satisfies $k=3, P=1$ and $L=2$. Later, Adamczewski and Bugeaud proved the Cobham conjecture. Thus, they deduce
Theorem \ref{MM-AB} by combining the results of Morton and Mourant with the Cobham conjecture.
\newline
\quad
The purpose of this paper is to extend Theorem 1 to uncountably many real numbers by the generalization of the sequence $(e_P^L(n))_{n\ge 0}$.\\
\quad
Now we give the following generalization of ${(e^L_P(n))}_{n\ge0}$. Let $d$ be a positive integer.
For a pattern $P:=p_1 \cdots p_d$ with $P\in {\{0,1,\cdots k-1\}}^{*}/{\{0\}}^{*}$, set the map $\mu_{P}: \mathbb{N} \longrightarrow \{0,1, \dots, L-1 \}$.
For the pattern $P$ and any non-negative integer $y$, we define the digital pattern counting function $d_1(n; p_d k^y+\cdots +p_1k^{y+d-1})$ as
\begin{displaymath}
d_1(n; p_d k^y+\cdots +p_1k^{y+d-1}):=\begin{cases} \mu_P(y) & 
 \mbox{  $p_d k^y+\cdots +p_1k^{y+d-1}$ is}\\ & 
\mbox{appeared in the base-$k$ representation of $n$} \\
                      0 & \mbox{Otherwise.}
                                             \end{cases}
\end{displaymath}
We define the sequence $(a(n))_{n\ge0}$ as
\begin{align}\label{LPRSRS}
    a(n) \equiv \sum_{y=0}^\infty d_1(n;p_d k^y+\cdots +p_1k^{y+d-1}) \pmod{L},
\end{align}
where $0\le a(n)\le L-1$ and $a(0)=0$. (For any non-negative integer $n$, the number of non-zero digit of $n$ is finite. Therefore, the value of right-hand sides \eqref{LPRSRS} is finite.) We call $(a(n))_{n\ge0}$ a generalized Rudin-Shapiro sequence of type $(L,k,P,\mu_{P})$, abbreviated as the $(L,k,P,\mu_{P})$-RS sequence.
In this paper, we generalize Theorem \ref{MM-AB}  as follows.
\begin{theorem}\label{main}
{\rm
Let $(a(n))_{n\ge0}$ be an $(L,k,P,\mu_{P})$-RS sequence with $2\le d$. Let $\beta$ be an integer with $\beta \ge L$. If there is no integer $A$ such that
\begin{align*}
\mu_{P}(y)=0     
\end{align*}
 for any integer $y$ with $ A\le y$, then $\sum_{n=0}^\infty \frac{a(N+nl)}{ {\beta}^{n+1}}$ $($ for all $N \ge 0$ and for all $l >0$ $)$ is a transcendental number $($ Moreover, an $S$-, $T$-, or  $Liouville$ number $)$. }
\end{theorem}
By uncountable many possibilities for the choice of $\mu_P$, we get uncountably many  transcendental numbers, which are  generated by the $(L,k,P,\mu_{P})$-RS sequences. Moreover, if $\sum_{n=0}^\infty \frac{a(n)}{ {\beta}^{n+1}}$ is a transcendental number,  then $\sum_{n=0}^\infty \frac{a(N+nl)}{ {\beta}^{n+1}}$ $($ for all $N \ge 0$ and for all $l >0$ $)$ is also a transcendental number. The proof of Theorem \ref{main} does not rest on pure periodicity of the periodic generalized Rudin-Shapiro sequence $(a(n))_{n\ge0}$ and the Cobham conjecture. Here we study non-periodicity of the subsequence $(a(N+nl))_{n\ge0}$ $($ for all $N \ge 0$ and for all $l >0$ $)$ of a generalized Rudin-Shapiro sequence $(a(n))_{n\ge0}$. Almost no generalized Rudin-Shapiro sequence $(a(n))_{n\ge0}$ is $k$-automatic (See Remark \ref{rem3}). Therefore, the proof of Theorem \ref{main} is different from the proof of Theorem \ref{MM-AB}. For the proof of Theorem \ref{main}, we introduce the notion of generalized $k$-regular sequences which is introduced as the generalization of $k$-automatic sequences (See Section $3$ and Appendices). Then we prove Theorem \ref{main} by combining Theorem \ref{period} (See Section $2$) with the combinatorial quantitative transcendence criterion established by Adamczewski and Bugeaud \cite{AB2}.
\begin{remark}{\rm(\cite{Mi1})
For the case $d=1$, we deal with in \cite{Mi1}. See Theorem $2$ in \cite{Mi1}.}
\end{remark}
\begin{remark}{\rm(\cite{Bu}) Here we introduce the definitions of Mahler's $S$-, $T$-, $U$-, and  $Liouville$ numbers. 
Let $\zeta$ be a real number. Let $w_n(\zeta)$ be the supremum of the real numbers $w$ for
which there exist infinitely many integer polynomials $R(x)=c_0+c_1 x+\cdots +c_n x^n$ of degree at most $n$ satisfying
$$
0<|R(\zeta)|\le \frac{1}{{H(R)}^w}
$$
where $H(R)=\max\{|c_i|\; | \; 0\le i\le n\}$.
We define the value $w(\zeta)$ as $w(\zeta):=\lim_{n=1}^\infty \sup \frac{w_n(\zeta)}{n}$.
We say that $\zeta$ is an 
\newline
\newline
$A$-number, if  $w(\zeta)=0$,
\newline
\newline
 $S$-number, if  $0<w(\zeta)<+\infty$,
\newline
\newline
 $T$-number, if  $w(\zeta)=+\infty$ and $w_n(\zeta)<+\infty$ for any positive integer $n$,
\newline
\newline
 $U$-number, if  $w(\zeta)=+\infty$ and $w_n(\zeta)=+\infty$ for some positive integer $n$.
\newline
\newline
Especially, we say that $\zeta$ is a $Liouville$ number, if  $w(\zeta)=+\infty$ and $w_1(\zeta)=+\infty$.}
\end{remark}
\begin{remark}\label{rem3}{\rm Now we show that almost no generalized Rudin-Shapiro sequence is $k$-automatic.
More precisely, we prove that if the map $\mu_{P}$ is non-periodic, then the $(L,k,P,\mu_{P})$-RS sequence $(a(n))_{n\ge0}$ is not $k$-automatic.
Assume that $(a(n))_{n\ge0}$ is $k$-automatic. By Theorem $6.6.2$ in \cite{AlS2}, there exists an integer $M$ such that  $$\{ (a(k^e n+j))_{n=0}^\infty\; | \;0\le e, \; 0\le j\le k^e-1  \}=\{ (a(k^e n+j))_{n=0}^\infty\; | 0\le e\le M, \; 0\le j\le k^e-1  \}.$$
Thus, there exist integers $e$ with $0\le e\le M$ and $j$ with $0\le j\le k^e-1$ such that,
\begin{align}\label{eq1.3}
a(k^{M+1} n)=a(k^e n+j).
\end{align}
Substitute $0$ for $n$ in \eqref{eq1.3},
we have $a(j)=0$. Let $y$ be a non-negative integer. Substitute $p_d k^{y+d+1}+\cdots +p_1k^{y+2d-1}$ for $n$ in \eqref{eq1.3}, we have
\begin{align}\label{ll}
&\mu_{P}(M+2+d+y)=a(k^{M+1}(p_d k^{y+d+1}+\cdots +p_1k^{y+2d}))=a(k^e (p_d k^{y+d+1}+\cdots +p_1k^{y+2d})+j)\end{align}
By the definition of $k^e (p_d k^{y+d+1}+\cdots +p_1k^{y+2d})+j$, the base $k$-representation of $k^e (p_d k^{y+d+1}+\cdots +p_1k^{y+2d})+j$ is as follows
\begin{align}\label{llpp}
 k^e (p_d k^{y+d+1}+\cdots +p_1k^{y+2d})+j=p_1 p_{2}\cdots p_d \underbrace {0 \cdots 0}_{d+y+1}\underbrace{\cdots }_{e}.
\end{align}
By \eqref{llpp} and $a(j)=0$, we have
\begin{align}\label{lll}
a(k^e (p_d k^{y+d+1}+\cdots +p_1k^{y+2d})+j)\equiv a(k^e (p_d k^{y+d+1}+\cdots +p_1k^{y+2d}))+a(j)=\mu_{P}(e+1+d+y).
\end{align}
From \eqref{ll} and \eqref{lll}, the map $\mu_{P}$ has period $M-e+1$.
}
\end{remark}
\quad
This paper is organized as follows. In Section $2$, we introduce the basic concepts of the periodicity of sequences, the lemma on the base $k$-representation of non-negative integers. We use this lemma on the base $k$-representation of non-negative integers and  the properties of generalized Rudin-Shapiro sequences to prove a necessary and sufficient condition for their periodicity. Furthermore, we prove that if the generalized Rudin-Shapiro sequence is not periodic then no subsequence $(a(N+nl))_{n\ge0}$ ( for all $N \ge 0$ and for all $l >0$) of the generalized Rudin-Shapiro sequences is periodic. (See Theorem \ref{period}.) 
In Section $3$, we introduce the notion of generalized $k$-regular sequences and we explain certain fractal structure of generalized $k$-regular sequences by using Theorem \ref{Cartier} in Appendix A for the proof of Theorem \ref{main}.
In Section $4$, first, we introduce the combinatorial quantitative transcendence criterion established by Adamczewski and Bugeaud \cite{AB2}. Then, by using certain fractal structure of generalized $k$-regular sequences, we prove that all generalized $k$-regular sequences, which take their values in finite sets, have sublinear complexity.
Next, we show that the generalized Rudin-Shapiro sequence can be regarded as a generalized $k$-regular sequence by using Theorem \ref{pattern} in Appendix B.
Finally, by applying the combinatorial quantitative transcendence criterion to the generalized non-periodic Rudin-Shapiro sequence $(a(n))_{n\ge0}$, which takes its values in $\{0,1,\ldots,\beta-1\}$, we show that $\sum_{n=0}^\infty a(N+nl){\beta}^{-n-1}$ is transcendental.
In Appendix A, we prove the fundamental properties of generalized $k$-regular sequences. We also prove properties that are not necessarily used in the proof of Theorem \ref{main}, because we think generalized $k$-regular sequences itself an interesting object. (Theorems \ref{Cartier} and \ref{fund}-$4$ are only used to prove Theorem \ref{main}.)
In Appendix B, we give the examples of generalized $k$-regular sequences related to the base-$k$ representation of integers and infinite chains of Mahler functional equations. (Theorem \ref{pattern} are only used to prove Theorem \ref{main}.) 
From Appendices, we think that the notion of generalized $k$-regular sequences is a centaral object of infinite chains of Mahler functional equations (For instance, see \cite{LP}, \cite{DuN}, \cite{Ta}, \cite{AmV}.
See also \eqref{8}, \eqref{sum} and \eqref{product}.) in the future. We also think the development of infinite chains of Mahler functional equations may give more the concrete construction method of uncountable many transcendental numbers than this paper.
\section{Non-periodicity of $(L,k,P,\mu_{P})$-RS sequences}\label{s:3}
\quad
First we introduce some definition and lemma for the non-periodicity of $(L,k,P,\mu_{P})$-RS sequences;
Let $(a(n))_{n\ge0}$ be a sequence with values in $\mathbb{C}$. $(a(n))_{n\ge0}$ is called {\it ultimately \; periodic} if there exist non-negative integers $N$ and $l>0$ such that 
\begin{align}
    a(n)=a(n+l)  \qquad (\;for \;all\; n \ge N). 
\end{align}
An {\it arithmetical\;  subsequence} of $(a(n))_{n\ge0}$ is $(a(N+tl))_{t\ge0}$, where $N \ge 0$ and  $l >0$. 
\begin{definition}
{\rm Let $(a(n))_{n\ge0}$ be a sequence with values in $\mathbb{C}$. The sequence $(a(n))_{n\ge0}$ is called  everywhere non-periodic if all of its arithmetical subsequences are non-constants.}
\end{definition}

\begin{lemma}{\rm(\cite{Mi1})}\label{ari}
{\rm
$(a(n))_{n\ge0}$ is everywhere non periodic 
if and only if no arithmetical subsequence of $(a(n))_{n\ge0}$ is ultimately periodic.}
\end{lemma}
The proof of this lemma is found in \cite{Mi1} (see Corollary $1$ in \cite{Mi1}).
\quad
The next lemma is the key Lemma to prove the non-periodicity of the arithmetical subsequences of an $(L,k,P,\mu_{P})$-RS sequence.
\begin{lemma} \label{multiple}{\rm(\cite{Mi1})}
{\rm
Notation is the same as for \eqref{eq1.1}.
For any non-negative integer $n$, we define the integer $q_m(n)$ as $q_m(n):=\min \{ q \;| s_{n,q}\neq 0\}$.
If $l$ and $t$ denote non-negative integers, then there exists an integer $x$ such that
\begin{align}
s_{xl,q_m(xl)}=1\;and\;  s_{xl,q}=0
\end{align}
for any integer $q$ with $q_m(xl)<q<q_m(xl)+t$.
\newline
\quad Furthermore, if $t'$ is another non-negative integer, then there exists an integer $X$ such that
\begin{align}
q_m(xl)=q_m(Xl), s_{Xl,q_m(Xl)}=1, \; and\;  s_{Xl,q}=0
\end{align}
for any  integer $q$ with $q_m(Xl)<q<q_m(Xl)+t'$.}
\end{lemma}
The proof of this lemma is found in \cite{Mi1} (see Lemma $4$ in \cite{Mi1}).
\newline
\quad
Now we prove the following theorem.
\begin{theorem}\label{period}
{\rm
Let $(a(n))_{n\ge0}$ be an $(L,k,P,\mu_{P})$-RS sequence with $2\le d$. The sequence $(a(n))_{n\ge0}$ is ultimately periodic if and only if there exists an integer A such that
\begin{align}\label{periodcon}
\mu_{P}(y)=0
\end{align}
for any integer $y$ with $ A\le y$.
\newline
\quad
Moreover, if the $(L,k,P,\mu_{P})$-RS sequence is not ultimately periodic, then no arithmetical
subsequence of it is ultimately periodic.}
\end{theorem}
\begin{proof}
{\rm
First, we show that if there exists an arithmetical
subsequence of $(a(n))_{n\ge0}$ which is ultimately periodic, then $\mu_P$ is ultimately zero, i.e. satisfies \eqref{periodcon}.
For any non-negative integer $n$, we define the integer $q_M(n)$ as $q_M(n):=\max \{ q \;| s_{n,q}\neq 0\}$.
Assume that there exists an arithmetical
subsequence of $(a(n))_{n\ge0}$ which is ultimately periodic and
\begin{align}\label{perip}
\# \{ y \;|\; \mu_{P}(y)\neq 0 \}=\infty.
\end{align}
From Lemma \ref{ari}, there exist non-negative integers $N$ and $l(0<l )$ such that 
\begin{align}\label{ariperiod}
a(N)=a (N+n l)  \qquad   (\forall n \in \mathbb{N}).
 \end{align}
For any integers $t$ with $t >0$ and $n$ with $n\ge q_M(N)+d+1$, the base $k$-representation of $N+k^n t l$ is as follows
\begin{align}\label{llpp1}
N+k^n t l=\cdots  \underbrace {0 \cdots 0}_{d}\underbrace{\cdots }_{q_M(N)+1}.
\end{align}
By \eqref{llpp1}, we get
 \begin{align}\label{keyeqp}
a (N+k^n t l) \equiv a (N)+a(k^n t l)   \; \; \;  \; {\rm for \; all}\;  n\ge q_M(N)+d+1, \;{\rm for \; all}\; t >0.
\end{align}
By \eqref{ariperiod} and \eqref{keyeqp}, we have
\begin{align}\label{keyeq}
a (k^n t l) =0   \; \; \;  \; {\rm for \; all}\;  n\ge q_M(N)+d+1, \;{\rm for \; all}\; t >0.
\end{align}
By Lemma \ref{multiple} and the assumption of $\mu_{P}(y)$, there exists an integer $x$ such that
\begin{align}\label{x}
\mu_{P}(q_m(xl))\neq 0,\;q_M(N)+d+1 <q_m(xl), \; s_{xl,q_m(xl)}=1 \;{\rm and}\; s_{xl,q}=0,  
\end{align}
for any integer $q$ with $q_m(xl)<q<q_m(xl)+3d+1$. For the pattern $P:=p_1 \cdots p_d$ with $P\in {\{0,1,\cdots k-1\}}^{*}/{\{0\}}^{*}$, let $R:=\sum_{j=0}^{d-1} p_{j+1} k^{d-1-j}$.
Moreover, by Lemma \ref{multiple}, there exists an integer $X$ such that
\begin{align}\label{X}
q_m(xl)=q_m(Xl), s_{Xl,q_m(Xl)}=1, \; and\;  s_{Xl,q}=0,  
\end{align}
for any integer $q$ with $q_m(Xl)<q<q_m(Xl)+2d+1+q_M(Rxl)$.\\
\quad
We consider the following two cases $1$ and $2$.
\newline
\newline
{\bf{Case}} $1$. $p_1\neq 0$. 
\newline
Let $U:=p_1 k^{d-1}$ and $W:=R-U$. By \eqref{x} and \eqref{X}, the base $k$-representations of $U x l$ and $WXl$ are as follows
\begin{align}
                       & Uxl= \cdots \underbrace {0 \cdots 0}_{2d+1}p_1  \underbrace {0 \cdots 0}_{q_m(xl)+d-1},
\\&   WXl= \cdots \underbrace {0 \cdots 0}_{2d+q_M(Rxl)+1}  p_{2} \cdots p_d \underbrace{ 0 \cdots \cdots \cdots 0}_{q_m(xl)}.   
\end{align}
By the definitions of $U x l$ and $WXl$, the base $k$-representation of $(Ux+WX)l$ is as follows
\begin{align}
                        (Ux+W X)l= \cdots  \underbrace {0 \cdots 0}_{2d+1} p_1 p_{2}\cdots p_d \underbrace {0 \cdots 0}_{q_m(xl)}.
\end{align}
By \eqref{keyeq}, the definitions of $Ux l$, $WXl$ and $(Ux+W X)l$, we get
\begin{align}
&a (U x l)=0,
\\&a (W X l)=0,
\\&a ((Ux+W X) l)=0.                     
\end{align}
From $(2.12)$-$(2.17)$ and $\mu_{P}(q_m(xl))\neq 0$, we have  $P = p_10\cdots0$.
\newline
\quad
Let $T_1:=p_1 k^d$. Then the base $k$-representations of $R x l$ and $T_1Xl$ are as follows
\begin{align}
                       & \; \; \; \; \;\;\;\;\;\; Rxl= \cdots  \underbrace {0 \cdots 0}_{2 d+1}p_1\underbrace {0 \cdots 0}_{q_m(xl)+d-1},
\\&   T_1Xl= \cdots \cdots\underbrace {0 \cdots 0}_{d+q_M(Rxl) +1}p_1 \underbrace {0 \cdots 0}_{q_m(xl)+d}.   
\end{align}
By the definitions of $R x l$ and $T_1Xl$, the base $k$-representation of $(Rx+T_1 X)l$ is as follows
\begin{align}
                        (Rx+T_1 X)l= \cdots \underbrace {0 \cdots 0}_{2 d} p_1 p_1 \underbrace {0 \cdots 0}_{q_m(xl)+d-1}.
\end{align}
From \eqref{keyeq} and the definition of $T_1 Xl$, we get
\begin{align}
   &a(T_1 X l)=0,\\
   &a((Rx+T_1 X) l)=0.
\end{align}
On the other hands, by $(2.18)$-$(2.22)$ and $p_1\neq 0$, we have
\begin{align}
   \mu_{P}(q_m(xl)+1)=0.                     
\end{align}
For any positive integer $s$, let $T_s:=p_1 k^{d+s-1}$.
From the base $k$-representations of $T_sxl$, $T_{s+1}Xl$ and the same process of $(2.18)$-$(2.22)$, we have
\begin{align}
   \mu_{P}(q_m(xl)+s)=0.                     
\end{align}
for any integer $s>0$.
This contradicts the fact $\# \{ y \;|\; \mu_{P}(y)\neq 0 \}=\infty$.
\newline 
\newline
{\bf{Case}} $2$.  $p_1= 0$. 
\newline
Let $J:=\min\{j | p_j\neq 0\}$. We put $Q:=\sum_{i=d-J}^{2d-1-J} k^i$. Then the base $k$-representations of $R x l$ and $QXl$ are as follows
\begin{align}
                       & \; \; \; \; \;\;\;\;\;\; Rxl= \cdots  \underbrace {0 \cdots 0}_{2 d+1}p_J \cdots p_d\underbrace{0 \cdots 0}_{q_m(xl)},
\\&   QXl= \cdots \cdots\underbrace {0 \cdots 0}_ {q_M(Rxl)+d}\underbrace {1 \cdots 1}_{d} \underbrace {0 \cdots 0}_{q_m(xl)+d-J+1}.   
\end{align}
By the definitions of $R x l$ and $QXl$, the base $k$-representation of $(Rx+Q X)l$ is as follows
\begin{align}
                        (Rx+QX)l=\cdots \cdots\underbrace {0 \cdots 0}_ {d}\underbrace {1 \cdots 1}_{d}p_J \cdots p_d\underbrace{0 \cdots 0}_{q_m(xl)}
\end{align}
By \eqref{keyeq}, the definitions of $Rx l$, $QXl$ and $(Rx+Q X)l$, we get
\begin{align}
&a (Rx l)=0,
\\&a(Q X l)=0,
\\&a ((Rx+Q X) l)=0.                     
\end{align}
By $(2.25)$-$(2.30)$ and $p_1= 0$, we have $\mu_{P}(q_m(xl))=0$.
This contradicts the fact $\mu_{P}(q_m(xl))\neq0$.
We proved that condition \eqref{periodcon} is necessary.
\newline
\quad
Finally, we show that \eqref{periodcon} gives the sufficient condition.
Assume that $\mu_P$ satisfies the condition \eqref{periodcon}. We define the integer $Y$ by $Y:=\max\{ y \;|\; \mu_{P}(y)\neq 0 \}$.
By the definition of $(a(n))_{n=0}^\infty$ (Note the case $p_1=0$), we have
\begin{align}
a(n)\equiv a(n)+a(tk^{Y+1+d})\equiv a(n+tk^{Y+1+d})
\end{align}
for  any integers $t>0$ and $n$ with $0\le n\le k^{Y}-1$.
Therefore, the sequence $(a(n))_{n\ge0}$ is ultimately periodic.
This completes the proof of Theorem \ref{period}.\hfill
$\square$
}
\end{proof}
\section{The certain fractal structure of generalized $k$-regular sequences}\label{s:1}
\quad In this section, we introduce the notion of of generalized $k$-regular sequences and we explain certain fractal structure of generalized $k$-regular sequences for the proof of Theorem \ref{main}.\\
\quad
First, we introduce the notion of generalized $k$-regular sequences. Let ${\bf{b}}:={(b(n ))}_{n\ge0}$ be a given sequence which takes its values in $\mathbb{C}$.
The sequence ${(b(n ))}_{n\ge0}$ is $k$-automatic if and only if the set $S$ defined by
\begin{align*} S:=\{{(b(k^en+j ))}_{n\ge0}\; | \; e\ge 0 , 0\le j \le k^e-1 \}\end{align*} 
is finite. Motivated by this result, in \cite{AlS1}, Allouche and Shallit introduced the notion of $k$-regular sequence generalizing that of $k$-automatic sequence as follows. A sequence ${(b(n ))}_{n\ge0}$ is
defined to be $k$-regular if the set $S$ is contained in $\mathbb{C}$-subspace generated by certain finite sequences.\\
\quad
Now, we introduce the notion of generalized $k$-regular sequences.
Let ${\bf{b}}:={(b(n ))}_{n\ge0}$ be a given sequence which takes its values in $\mathbb{C}$.
For any non-negative integer $e$, set $$S_e({\bf{b}}):=\{{(b(k^en+j ))}_{n\ge0}\; |\;0\le j \le k^e-1 \}.$$
\begin{definition}\label{defofk-pro}{}
{\rm A sequence ${(b(n ))}_{n\ge0}$ is {\it generalized $k$-regular} if and only if there exists an integer $d$ such that
the dimension of the $\mathbb{C}$-subspace generated by $S_e({\bf b})$ is at most $d$ for all $e\ge 0$.
To show the role of $d$ more precisely, generalized $k$-regular is also called
{\it generalized $(k,d)$-regular}.}
\end{definition}
\quad
Finally, we seek certain fractal structure of generalized $k$-regular sequences by using Theorem \ref{Cartier} in Appendix A.
Notation is the same as for \eqref{eq1.1} and Theorem \ref{Cartier}.
For any generalized ($k$, $d$)-regular sequence ${\bf{b}}$, non-negative integers $e$ and $n$, we define the vector letter ${\bf{B}}_{e}(n)$ and the infinite vector word ${\bf{B}}_{e}$ as follows
\begin{align}
 & {\bf{B}}_{e}(n):=\left(
    \begin{array}{c}
       b_{e,1}(n)\\
       {b_{e,2}(n)}\\
      \vdots \\
      {b_{e,d}(n)}
    \end{array}
  \right),\\
&{\bf{B}}_{e}:=\left(
    \begin{array}{c}
       b_{e,1}(0)b_{e,1}(1)\cdots b_{e,1}(n)\cdots\\
       b_{e,2}(0)b_{e,2}(1)\cdots b_{e,2}(n)\cdots\\
      \vdots \\
      b_{e,d}(0)b_{e,d}(1)\cdots b_{e,d}(n)\cdots
    \end{array}
  \right)={\bf{B}}_{e}(0)\cdots{\bf{B}}_{e}(n)\cdots.
  \end{align}
Especially, the generalized ($k$, $d$)-regular sequence ${\bf{b}}$ is the first coordinate of the infinite vector word ${\bf{B}}_{0}$. \\
\quad
For any positive integers $e$ and $m$ with $0\le m\le k^e-1$, we define the $d\times d$ matrix $D_{m,e}$ as
\begin{align}\label{matrix}
D_{m,e}:=C_{s_{m,0},0}C_{s_{m,1},1}\cdots C_{s_{m,e-1},e-1}.
\end{align}
By Theorem \ref{Cartier}, for any positive integer $e$, we have
\begin{align}\label{HFi}
&{\bf{B}}_{e-1}=
\underbrace{C_{0,e-1}{\bf{B}}_{e}(0) C_{1,e-1}{\bf{B}}_{e}(0)\cdots C_{k-1,e}{\bf{B}}_{e}(0)}_{k}
\underbrace{C_{0,e-1}{\bf{B}}_{e}(1) C_{1,e-1}{\bf{B}}_{e}(1)\cdots C_{k-1,e}{\bf{B}}_{e}(1)}_{k}\nonumber \\
&\cdots \underbrace{C_{0,e-1} {\bf{B}}_{e}(n) C_{1,e-1}{\bf{B}}_{e}(n)\cdots C_{k-1,e}{\bf{B}}_{e}(n)}_{k}\cdots\cdots \cdots \cdots \cdots .\end{align}
Moreover, by using Theorem \ref{Cartier} repeatedly and \eqref{matrix}, for any positive integer $e$, we have
\begin{align}\label{HF}
&{\bf{B}}_{0}=
\underbrace{D_{0,e}{\bf{B}}_{e}(0) D_{1,e}{\bf{B}}_{e}(0)\cdots D_{k^e-1,e}{\bf{B}}_{e}(0)}_{k^e}
\underbrace{D_{0,e}{\bf{B}}_{e}(1) D_{1,e}{\bf{B}}_{e}(1)\cdots D_{k^e-1,e}{\bf{B}}_{e}(1)}_{k^e}\nonumber \\
&\cdots \underbrace{D_{0,e} {\bf{B}}_{e}(n) D_{1,e} {\bf{B}}_{e}(n)\cdots D_{k^e-1,e} {\bf{B}}_{e}(n)}_{k^e}\cdots\cdots \cdots \cdots \cdots .\end{align}
This representation \eqref{HF} of ${\bf{B}}_{0}$ is the fractal structure of the generalized ($k$, $d$)-regular sequence $\bf{b}$.
This fractal structure is the key property for applying the combinatorial quantitative transcendence criterion (See Theorem \ref{cqtc} in Section $4$.) to the generalized $k$-regular sequences, which take their values in finite set.
\section{Proof of Theorem \ref{main}}\label{s:1}
\quad In this section, we prove Theorem \ref{main} by using the results in previous sections.
First, we introduce the combinatorial quantitative transcendence criterion established by Adamczewski and Bugeaud \cite{AB2}.
Let $\beta$ be an integer greater than $1$ and ${\bf{a}}:={(a(n))}_{n\ge0 }$ be a sequence on $\{ 0,1,\cdots, \beta-1\}$. For any non-negative integer $m$ and the sequence ${(a(n))}_{n\ge0 }$, we define the complexity function $p_{\bf{a}}(m)$ as
$$
p_{\bf{a}}(m)=\#\{ a(j)a(j+1)\cdots a(j+m-1) \; | j\ge0\}.
$$
We define the subset $\mathfrak{CL}$ of real numbers as
\begin{align*}
\mathfrak{CL}:=\{ \zeta \in \mathbb{R}\; | \;\mbox {there exists a base $b$ such that} \;p_{\bf{a}}(m)=O (m) \;\mbox {where} \;\zeta=\sum_{n=0}^\infty \frac{a(n)}{ {\beta}^{n+1}}  \}.
\end{align*}
Adamczewski and Bugeaud \cite{AB2} proved the following theorem by using quantitative Subspace Theorem.
\begin{theorem}{\rm(\cite{AB2})}\label{cqtc}
{\rm
Let $\zeta$ be a real number in $\mathfrak{CL}$. Then $\zeta$ is a rational, $S$-, $T$-, or  $Liouville$ number.}
\end{theorem}
\quad
Now we prove that all generalized $k$-regular sequences, which take their values in finite sets, have sublinear complexity. (i.e., such that $p_{\bf{a}}(m)=O(m)$). Theorem \ref{geCob} is a generalization of Theorem $10.3.1$ in \cite{AlS2}. Theorem \ref{geCob} is proved in a similar method of the proof of Theorem $10.3.1$ in \cite{AlS2}.
\begin{theorem}\label{geCob}
{\rm
Let $\beta$ be an integer greater than $1$ and ${\bf{a}}:={(a(n))}_{n\ge0}$ be a sequence which takes its values on $\{0,1,\cdots,\beta-1\}$.
Moreover, there exists an injective coding $\tau:{\{0,1,\cdots,\beta-1\}}^{*}\rightarrow {\mathbb{C}}^{*}$ such that ${\bf{b}}:={(b(n))}_{n\ge0}:={(\tau(a(n)))}_{n\ge0}$ is a generalized ($k$, $d$)-regular sequence which takes its values in ${\mathbb{C}}$.
Then, for any non-negative integer $m$,
\begin{align}\label{sublinear}
p_{\bf{a}}(m)\le  {(\beta+1)}^{2d} km.
\end{align}}
\end{theorem}
\begin{proof}
{\rm
Notation is the same as for Section $3$ and Theorem \ref{Cartier}.
By Theorem \ref{Cartier}, we assume that any sequence ${\bf{b}}_{e,l}$ ($e\ge0$, $1\le l \le d$) is in $S_e({\bf{b}})$ and zero sequence. Therefore, ${\bf{b}}_{e,l}$ ($e\ge0$, $1\le l \le d$) takes its values on at most $\beta+1$ values.
Let $m$ be a non-negative integer, $c_1\cdots c_m$ be a $m$ length word in ${\bf{a}}$ and $y$ be an integer with $k^{y-1}\le m < k^{y}$.
From the property \eqref{HF} of generalized $k$-regular sequences, the word $\tau(c_1\cdots c_m)$ is included in the first coordinate of the following vector word
\begin{align}\label{comp}
\underbrace{D_{0,y}{\bf{B}}_{y}(n) D_{1,y}{\bf{B}}_{y}(n)\cdots D_{k^y-1,y}{\bf{B}}_{y}(n)}_{k^y}
\underbrace{D_{0,y}{\bf{B}}_{y}(n+1) D_{1,y}{\bf{B}}_{y}(n+1)\cdots D_{k^y-1,y}{\bf{B}}_{y}(n+1)}_{k^y}
\end{align}
where $n$ is some positive integer.
By the injectivity of $\tau$ and \eqref{comp}, we get
\begin{align}
p_{\bf{a}}(m)=p_{\bf{b}}(m)\le  {(\beta+1)}^{2d}   k^y \le  {(\beta+1)}^{2d} km.
\end{align}
\hfill
$\square$}
\end{proof}
From Theorem \ref{cqtc} and Theorem \ref{geCob}, we generalizes Cobham conjecture, which was proved by Adamczewski and Bugeaud (See Theorem $2$ in \cite{AB1}.), as follows.
\begin{corollary}\label{geCob2}
{\rm
Let $\beta$ be an integer greater than $1$ and ${\bf{a}}:={(a(n))}_{n\ge0}$ be a sequence which takes its values on $\{0,1,\cdots,\beta-1\}$.
Moreover, there exists an injective coding $\tau:{\{0,1,\cdots,\beta-1\}}^{*}\rightarrow {\mathbb{C}}^{*}$ such that ${\bf{b}}:={(b(n))}_{n\ge0}:={(\tau(a(n)))}_{n\ge0}$ is a generalized ($k$, $d$)-regular sequence which takes its values in ${\mathbb{C}}$.
Then $ \sum_{n=0}^\infty \frac{a(n)}{\beta^{n+1}}$ is a rational, $S$-, $T$-, or  $Liouville$ number.}
\end{corollary}
Now, we explain that the $(L,k,P,\mu_{P})$-RS sequence can be regarded as a generalized $k$-regular sequence.
Notation is the same as for Theorem \ref{pattern} in Appendix B.
For any integer $j$ with $0\le \; j \le L-1$, we define the injective coding $\tau:{\{0,1,\cdots,L-1\}}^{*}\rightarrow \mathbb{C}^{*}$ as $\tau(j)=\exp \frac{2\pi \sqrt{-1} j}{L}$. (For the definition of coding, see $9$p in \cite{AlS2}.)
For the $(L,k,P,\mu_{P})$-RS sequence $(a(n))_{n\ge0}$ with $P=p_1\cdots p_d$, we define the sequence $(\tau(a(n)))_{n\ge0}$. This replaced sequence $(\tau(a(n)))_{n\ge0}$ is contained in the class of sequences defined by \eqref{mul} in Appendix B. Indeed, we define the map $\mu_{P}^{'}: \mathbb{N} \longrightarrow \mathbb{C}$ by $\mu_{P}^{'}=\tau \circ\mu_{P}^{}$.
By the definition of the $(L,k,P,\mu_{P})$-RS sequence $(a(n))_{n\ge0}$, $\mu_{P}^{'}=\tau \circ\mu_{P}^{}$ and $\tau(0)=1$, we have
\begin{align}
&\tau(a(n))=\tau(\sum_{y=0}^\infty d_1(n;p_d k^y+\cdots +p_1k^{y+d-1}))\nonumber \\
&=\prod_{y=0}^\infty \tau(d_1(n;p_d k^y+\cdots +p_1k^{y+d-1}))=\prod_{y=0}^\infty d_3(n;p_d k^y+\cdots +p_1k^{y+d-1}).
\end{align}
Then we get the following corollary.
\begin{Prop}\label{patternkpro}
{\rm
Let $(a(n))_{n\ge0}$ be an $(L,k,P,\mu_{P})$-RS sequence. Let $\tau$ be an injective coding $\tau:{\{0,1,\cdots,L-1\}}^{*}\rightarrow \mathbb{C}^{*}$ defined by $\tau(j)=\exp \frac{2\pi \sqrt{-1} j}{L}$. Then the sequence ${(\tau(a(n)))}_{n\ge0}$ is generalized $k$-regular.}
\end{Prop}\quad
Finally, we prove Theorem \ref{main}. 
\begin{proof}
{\rm From Proposition \ref{patternkpro} and Theorem \ref{fund}-4 in Appendix A, all arithmetical subsequences of the $(L,k,P,\mu_{P})$-RS sequence satisfy the assumption of Corollary \ref{geCob2}.
By combine Corollary \ref{geCob2} with Theorem \ref{period}, we deduce Theorem \ref{main}. \hfill
$\square$
}
\end{proof}
\section{Appendix A: the properties of generalized $k$-regular sequences}\label{s:1}
In this section, we gather fundamental results for
generalized $k$-regular sequences, which generalize corresponding results for $k$-regular sequences.\\
\quad
We here recall the Cartier operators. For any sequence
${(b(n))}_{n\ge0 }$, the Cartier operators ${\Delta}_j={\Delta}_{j,k}$ with $0\le j\le k-1$, act on $\bf{b}$ $:={(b(n))}_{n\ge0 }$ as 
$${\Delta}_j {\bf{b}}:={(b(kn+j))}_{n\ge0 }.$$
It can be easily seen that, for any non-negative integer $e$, we have
\begin{align}\label{s2}
\bigcup_{0\le j<k} 
\bigcup_{{\bf{c}}\in S_{e}({\bf{b}})} {\Delta}_j {\bf{c}}=S_{e+1}({\bf{b}}).
\end{align}
The definition of generalized ($k$, $d$)-regular sequences together with this relation implies the
following result, which partly generalizes Theorem $5.1.1$ in \cite{N}.
\begin{theorem}\label{Cartier}
{\rm A sequence ${\bf{b}}={(b(n))}_{n\ge0 }$ is generalized ($k$, $d$)-regular if and only if there exist sequences 
${{\bf{b}}_{e,u}}={(b_{e,u}(n))}_{n\ge0 }$ 
$($ $e\ge0$, $1\le u \le d$ $)$ with ${{\bf{b}}_{0,1}}={\bf{b}}$ and $d\times d$ matrices $C_{j,e}$ $($ $e\ge0$, $0\le j \le k-1$ $)$ whose entries are in $\mathbb{C}$ such that
\begin{align}\label{eqm}
{\Delta}_j\left(
    \begin{array}{c}
      {\bf{b}}_{e,1}\\
      {\bf{b}}_{e,2}\\
      \vdots \\
     {\bf{b}}_{e,d}
   \end{array}
  \right)=C_{j,e}
  \left(
    \begin{array}{c}
      {\bf{b}}_{e+1,1}\\
      {\bf{b}}_{e+1,2}\\
      \vdots \\
     {\bf{b}}_{e+1,d}
   \end{array}
  \right).
  \end{align}}
\end{theorem}
\begin{proof}
{\rm
Assume that  ${\bf{b}}$ is a generalized ($k$, $d$)-regular sequence.
For any non-negative integer $e$, 
we denote by ${\hat{S}}_e({\bf{b}})$ the $\mathbb{C}$-subspace generated by $S_e({\bf{b}})$. Then, by the definition of generalized ($k$, $d$)-regular sequences, the dimension $d_{e}$ of ${\hat{S}}_e({\bf{b}})$ is at most $d$ and ${\hat{S}}_e({\bf{b}})$ can be expressed as 
\begin{align}
{\hat{S}}_e({\bf{b}})=<{\bf{b}}_{e,1},{\bf{b}}_{e,2},\cdots{\bf{b}}_{e,d}>\;\; (e\ge0),
 \end{align} 
 where the right-hand side denote the $\mathbb{C}$-subspace generated by ${\bf{b}}_{e,u}\in S_e({\bf{b}})$ ($1\le u\le d_e$) with ${\bf{b}}_{0,1}=\bf{b}$ (We can take the largest linearly independent tuples of $S_e({\bf{b}})$.) and ${\bf{b}}_{e,u}=0$ ($d_e+1\le l\le d$ ). Then, by (\ref{s2}), we have ${\Delta}_j {\bf{b}}_{e,u}\in {{S}}_{e+1}({\bf{b}})\subset {\hat{S}}_{e+1}({\bf{b}})$
$($ $e\ge0$, $0\le j \le k-1$, $1\le u\le d$ $)$. Therefore, the desired relations \eqref{eqm} hold for these ${\bf{b}}_{e,u}$ with certain $C_{j,e}$.
\newline
\quad
Conversely, assume that \eqref{eqm} hold. Then the repeated application of \eqref{eqm} implies that 
\begin{align}
{{S}}_e({\bf{b}})\subset <{\bf{b}}_{e,1},{\bf{b}}_{e,2},\cdots{\bf{b}}_{e,d}>\;\; (e\ge0).
 \end{align} 
This completes the proof.
\hfill
$\square$}
\end{proof}
Here, for any sequences ${\bf{b}}={(a(n ))}_{n\ge0}$ and ${\bf{c}}={(b(n ))}_{n\ge0}$, the canonical convolution ${\bf{b}}*{\bf{c}}={((b*c)(n ))}_{n\ge0}$ of ${\bf{b}}$ and ${\bf{c}}$ is defined by
\begin{align*}
(b*c)(n )=\sum_{i=0}^n b(i)c(n-i)  \;\; (n\ge0).
\end{align*}
Moreover, for any sequence ${\bf{b}}={(a(n ))}_{n\ge0}$,  ${\bf{(b)}_{-1}}$ is defined by ${\bf{(b)}_{-1}}:={(b(n-1))}_{n\ge0}$ with $b(-1)=0$.
We prove the following closer properties of generalized $k$-regular sequences. Theorem \ref{fund} proved in a similar method of the proof of Theorems $2.5$, $2.6$ and $3.1$ in \cite{AlS1}.
\begin{theorem}\label{fund}
{\rm Let $\bf{b}$ $:={(b(n))}_{n\ge0 }$ be generalized ($k$, $d_1$)-regular and $\bf{c}$ $:={(c(n))}_{n\ge0 }$ be generalized ($k$, $d_2$)-regular.
Then the following sequences are also generalized $k$-regular sequences.
\newline
$1$. $\bf{b}$ $+$ $\bf{c}$$:={(b(n)+c(n))}_{n\ge0}$;
\newline
$2$. $\bf{b}\cdot \bf{c}$ $:={(b(n)\cdot c(n))}_{n\ge0}$;
\newline
$3$. $\alpha\bf{b}$ $:={(\alpha\cdot b(n))}_{n\ge0}$ for any $\alpha\in \mathbb{C}$;
\newline
$4$. ${\bf{b}}_{l,N}$$:={(b(N+ln))}_{n\ge0}$ for any integers $N$ and $ l$ with $N\ge0$ and $l\ge1$.}
\newline
$5$. ${\bf{b}}*{\bf{c}}$.
\end{theorem}
\begin{proof}
{\rm 
For any non-negative integer $e$, 
we denote by ${\hat{S}}_e({\bf{b}})$ and ${\hat{S}}_e({\bf{c}})$ the $\mathbb{C}$-subspace generated by $S_e({\bf{b}})$ and $S_e({\bf{c}})$. ${\hat{S}}_e({\bf{b}})$ and ${\hat{S}}_e({\bf{c}})$can be expressed as 
\begin{align}\label{psles}
&{\hat{S}}_e({\bf{b}})=<{\bf{b}}_{e,1},{\bf{b}}_{e,2},\cdots{\bf{b}}_{e,d_1}>\;\; (e\ge0),\nonumber\\&
{\hat{S}}_e({\bf{c}})=<{\bf{c}}_{e,1},{\bf{c}}_{e,2},\cdots{\bf{c}}_{e,d_2}>\;\; (e\ge0)
 \end{align}
where the right-hand sides denote the $\mathbb{C}$-subspace generated by ${\bf{b}}_{e,u}$ ($1\le u\le d_1$), ${\bf{c}}_{e,u}$ ($1\le u\le d_2$).\\
\quad
We prove the first statement. Let $\bf{d}$ $:={(a(n)+b(n))}_{n\ge0 }$. By \eqref{psles}, we have
\begin{align}\label{psles1}
{\hat{S}}_e({\bf{d}})=<{\bf{b}}_{e,1},{\bf{b}}_{e,2},\cdots{\bf{b}}_{e,d_1},{\bf{c}}_{e,1},{\bf{c}}_{e,2},\cdots{\bf{c}}_{e,d_2}>\;\; (e\ge0).
 \end{align}
Therefore, the first statement is proved.\\
\quad
We prove the second statement.
Let $\bf{f}$ $:={(a(n)\cdot b(n))}_{n\ge0 }$. By \eqref{psles}, we have
\begin{align}\label{psles1}
{\hat{S}}_e({\bf{f}})=<{\bf{b}}_{e,i}\cdot {\bf{c}}_{e,j} : 1\le i\le d_1, 1\le j\le d_2>\;\; (e\ge0).
 \end{align}
Therefore, the second statement is proved.\\
\quad
We prove the third statement. The third statement is trivial.\\
\quad
We prove the fourth statement.
For any non-negative integers $e$ and $j$ with $0\le j\le k^e-1$, set $N+jl=:k^e s +M$ where $0\le M\le k^e-1$. Then we have
\begin{align}\label{w1}
&{b}_{l,N}(k^en +j)=b(N+l(k^en +j) )=b(k^e(ln+s)+M).
\end{align}
From $0\le j\le k^e-1$, we have $0\le s\le N+l$.
For any non-negative integers $e$, $u$ and  $s$ with $1\le u \le d$ and $0\le s\le N+l$,
 we define ${\bf{b}}_{e,u,l,s}$ as ${\bf{b}}_{e,u,l,s}:={(b_{e,u}(ln+s))}_{n\ge0}$.
By Theorem \ref{Cartier}, we also have
\begin{align}\label{w2}
{(b(k^e(ln+s)+M))}_{n=0 }^\infty \in <{\bf{b}}_{e,1,l,s},\cdots,{\bf{b}}_{e,d,l,s}>.
\end{align}
By \eqref{w1} and \eqref{w2}, the set ${{S}}_{e}({\bf{b}}_{l,N})$ is contained in the $\mathbb{C}$-subspace generated by ${\bf{b}}_{e,u,l,s}$ ($1\le u\le d_1$, $0\le s\le N+l$).
The fourth statement is proved.
\\
\quad
Finally, we prove the fifth statement.
We define the $\mathbb{C}$-subspace ${\tilde{S}}_e({\bf{b}}*{\bf{c}})$ as 
\begin{align}\label{w2333}
{\tilde{S}}_e({\bf{b}}*{\bf{c}}) = <{\bf{b}}_{e,i}*{\bf{c}}_{e,j},({{\bf{b}}_{e,i}}* {{\bf{c}}_{e,j}})_{\bf{-1}}: 1\le i\le d_1, 1\le j\le d_2>\;\; (e\ge0).
\end{align}
For any integers $j_1$ and $j_2$ with $0\le j_1,j_2\le k-1$, let ${\bf{b}}^{e,j_1}:={(b^{e,j_1}(n))}_{n\ge0}={(b(k^en+j_1))}_{n\ge0}$ and ${\bf{c}}^{e,j_2}:={(c^{e,j_2}(n))}_{n\ge0}={(c(k^en+j_2))}_{n\ge0}$.
By \eqref{psles}, we get
\begin{align}\label{w2334}
{\bf{b}}^{e,j_1}*{\bf{c}}^{e,j_2},({\bf{b}}^{e,j_1}*{\bf{c}}^{e,j_2})_{\bf{-1}}\in {\tilde{S}}_e({\bf{b}}*{\bf{c}}).
\end{align}
For any integer $l$ with $0\le l\le k^e-1$, we have
\begin{align}\label{w2339}
&\sum_{0\le j_1\le l} ({{b}}^{e,j_1}*{{c}}^{e,l-j_1})(n)+\sum_{l< j_1<k^e} ({{b}}^{e,j_1}*{{c}}^{e,k^e+l-j_1})(n-1)=\nonumber\\&
\sum_{0\le j_1\le l} \sum_{i+j=n}b(k^e i+j_1)c(k^e j +l-j_1)+\sum_{l< j_1<k^e}\sum_{i+j=n-1}b(k^e i+j_1)c(k^e j+k^e+l-j_1)
=\nonumber\\&
\sum_{i+j=k^en+l} b(i)c(j)=
(b*c)(k^en +l).
\end{align}
From \eqref{w2333}-\eqref{w2339}, The fifth statement is proved.
\hfill
$\square$
}
\end{proof}
Theorem \ref{fund}-5 shows that the set of the generating series of generalized $k$-regular series forms a ring under the usual addition and product.\\
\quad
The following generalization of Theorem $5.1.2$ in \cite{N} characterizes generalized ($k$, $d$)-regular sequences by their generating series.
\begin{theorem}\label{th2}
{\rm A sequence ${(b(n ))}_{n\ge0}$ is generalized ($k$, $d$)-regular if and
only if, for any non-negative integer $e$, there exist $d$ power sereis $f_{e,1}(z),f_{e,2}(z),\cdots,f_{e,d}(z)\in K[[z]]$ with
$f_{0,1}(z)=\sum_{n=0}^\infty b(n)z^n$ and a $d\times d$ matrix $A_e(z)$ whose entries are polynomials
in $z$ of degrees less than $k$ with coefficients in $K$ such that
\begin{align}\label{8}
&\left(
    \begin{array}{c}
      f_{e,1}(z) \\
      f_{e,2}(z) \\
      \vdots \\
      f_{e,d}(z)
    \end{array}
  \right)=A_e(z)\left(
    \begin{array}{c}
      f_{e+1,1}(z^k) \\
      f_{e+1,2}(z^k)\\
      \vdots \\
      f_{e+1,d}(z^k)
   \end{array}
  \right) \;\;(e\ge0).
  \end{align} }
\end{theorem}
\begin{proof}
{\rm
Let us denote $C_{j,e}$ given in (\ref{eqm}) by $(c_{j,e}(p,q))_{1\le p,q \le d}$. Then we can see the varidity of (\ref{eqm}) is equivalent to that of (\ref{8}) with
\begin{align}
A_e(z)=(P_{e,p,q}(z))_{1\le p,q \le d}, \quad P_{e,p,q}(z)=\sum_{j=0}^{k-1}c_{j,e }(p,q)z^j \;\; (e\ge0).
\end{align}
The proof is the same as that of Theorem $5.1.2$ in \cite{N}.
\hfill
$\square$}
\end{proof}

\section{Appendix B: the examples of generalizaed $k$-regular seqeunces}\label{s:1}
In this section, we give several examples of generalized $k$-regular sequences.
\subsection{Generalized $k$-regular sequences related to the base $k$-representation of integers}
\quad In this subsection, we give examples of generalized $k$-regular sequences related to the base-$k$ representation of integers.\\
\quad
Gel'fond \cite{Gel} introduced the two functions related with the base $k$-representation as follows.
A sequence ${(a(n))}_{n\ge0}$ is $k$-additive if and only if, for any non-negative integers $e$, $n$ and $j$ with $0\le j\le k^e-1$, ${(a(n))}_{n\ge0}$ satisfies the following additive relation
\begin{align}\label{4}a(k^en+j )=a(k^en)+a(j)
\end{align}
and $a(0)=0$.
A sequence ${(a(n))}_{n\ge0}$ is $k$-multiplicative if and only if, for any non-negative integers $e$, $n$ and $j$ with $0\le j\le k^e-1$, ${(a(n))}_{n\ge0}$ satisfies the following multiplicative relation
\begin{align}\label{5}a(k^en+j )=a(k^en)a(j)
\end{align}
and $a(0)=1$.\\
\quad We here note that the set of $k$-regular sequences and the set of $k$-additive sequences (resp., $k$-multiplicative sequences) has large intersection. Indeed, for any non-negative integers $e$, $n$ and $j$ with $0\le j\le k^e-1$, if a $k$ additive sequence (resp.,  $k$-multiplicative sequence) ${(a(n ))}_{n\ge0}$ satisfies $a(k^en+j )=a(n)+a(j)$ (resp., $a(k^en+j )=a(n)a(j)$), then ${(a(n ))}_{n\ge0}$ is a $k$-regular sequence.
\\
\quad
In this appendix, we explain that the notion of generalized $k$-regular sequence unifes that of $k$-regular, $k$-additive and $k$-multiplicative sequence.\\
\quad
For a pattern $P:=p_1 \cdots p_d$ with $P\in {\{0,1,\cdots k-1\}}^{*}/{\{0\}}^{*}$, set the map $\mu_{P}^{'}: \mathbb{N} \longrightarrow \mathbb{C}$.
For the pattern $p_1 \cdots p_d$ and any non-negative integer $y$, we define the counting functions $d_2(n; p_d k^y+\cdots +p_1k^{y+d-1})$ and $d_3(n; p_d k^y+\cdots +p_1k^{y+d-1})$ as
\begin{displaymath}
d_2(n; p_d k^y+\cdots +p_1k^{y+d-1}):=\begin{cases} \mu_{P}^{'}(y)& \mbox{  $p_d k^y+\cdots +p_1k^{y+d-1}$ is}\\ & \mbox{appeared in the base-$k$ representation of $n$  } \\
                      0 & \mbox{Otherwise,}
                                             \end{cases}
\end{displaymath}
\begin{displaymath}
d_3(n; p_d k^y+\cdots +p_1k^{y+d-1}):=\begin{cases} \mu_{P}^{'}(y)& \mbox{  $p_d k^y+\cdots +p_1k^{y+d-1}$ is } \\ & \mbox{ appeared in the base-$k$ representation of $n$  } \\
                      1 & \mbox{Otherwise.}
                                             \end{cases}
\end{displaymath}
We define the sequence ${\bf{a}}=(a(n))_{n\ge0}$ as
\begin{align}\label{ad}
    a(n) = \sum_{y=0}^\infty d_2(n;p_d k^y+\cdots +p_1k^{y+d-1}).
\end{align}
We also define the sequence ${\bf{b}}=(b(n))_{n\ge0}$ as
\begin{align}\label{mul}
    b(n) = \prod_{y=0}^\infty d_3(n;p_d k^y+\cdots +p_1k^{y+d-1}).
\end{align}
For any non-negative integer $n$, the number of non-zero digit of $n$ is finite. Therefore, the values of right-hand sides \eqref{ad} and \eqref{mul} are finite.
\begin{theorem}\label{pattern1}
{\rm Let ${\bf{a}}$ be given in (\ref{ad}).
Then ${\bf{a}}$ is generalized $k$-regular.}
\end{theorem}
\begin{proof}
{\rm
By the definition of ${\bf{a}}$, for any non-negative integers $e$, $l$ and $j$ with $0\le l\le k^{d-1}-1$ and $0\le j\le k^{e}-1$, we have
\begin{align}\label{add}
&a(k^{e}(k^{d-1}n+l)+j)=\sum_{y=0}^\infty d_2(k^{e}(k^{d-1}n+l)+j;p_d k^y+\cdots +p_1k^{y+d-1})=\nonumber \\
&\sum_{y=e}^\infty d_2(k^{e}(k^{d-1}n+l)+j;p_d k^y+\cdots +p_1k^{y+d-1}) + \sum_{y=0}^{e-1} d_2(k^{e}(k^{d-1}n+l)+j;p_d k^y+\cdots +p_1k^{y+d-1}).\end{align}
By the definitions of $k^{e}(k^{d-1}n+l)+j$, the base $k$-representation of $k^{e}(k^{d-1}n+l)$ is as follows
\begin{align}\label{reprep}
           k^{e}(k^{d-1}n+l)+j= \underbrace { k^{e+d-1}n}_{}\underbrace {k^el}_{d-1} \underbrace {j}_{e}.
\end{align}
By the definition of $d_1$ and \eqref{reprep}, we have 
\begin{align}\label{add3}
\sum_{y=e}^\infty d_2(k^{e}(k^{d-1}n+l)+j;p_d k^y+\cdots +p_1k^{y+d-1})=\sum_{y=e}^\infty d_2(k^{e}(k^{d-1}n+l);p_d k^y+\cdots +p_1k^{y+d-1}).\end{align}
For any non-negative integers $e$ and $l$ with $0\le l\le k^{d-1}-1$, we define the sequence ${{\bf{a}}_{e,l}}={(a_{e,l}(n))}_{n\ge0 }$ as $$a_{e,l}(n)=\sum_{y=e}^\infty d_2(k^{e}(k^{d-1}n+l);p_d k^y+\cdots +p_1k^{y+d-1}).$$
By the definition of $d_2$ and \eqref{reprep}, we have
\begin{align}\label{add4}
\sum_{y=0}^{e-1} d_2(k^{e}(k^{d-1}n+l)+j;p_d k^y+\cdots +p_1k^{y+d-1})=\sum_{y=0}^{e-1} d_2(k^{e}l+j;p_d k^y+\cdots +p_1k^{y+d-1}).\end{align}
For any non-negative integers $e$, $l$ and $j$ with $0\le l\le k^{d-1}-1$ and $0\le j\le k^{e}-1$, we define the constant $c_{e,l,j}$ as $$c_{e,l,j}=\sum_{y=0}^{e-1} d_2(k^{e}l+j;p_d k^y+\cdots +p_1k^{y+d-1}).$$
By \eqref{add}, \eqref{add3}, \eqref{add4}, we get
\begin{align}\label{add5}
a(k^{e}(k^{d-1}n+l)+j)=c_{e,l,j}\times 1+a_{e,l}(n).
\end{align}
By \eqref{add5}, the set $S_{e+d-1}({\bf{b}})$ is contained in the $\mathbb{C}$-subspace generated by the $k^{d-1}+1$ sequences ${{\bf{a}}_{e,l}}$ ($0\le l\le k^{d-1}-1$) and ${(1)}_{n\ge0}$. Therefore, the sequence ${\bf{a}}$ is generalized $k$-regular.
\hfill
$\square$}
\end{proof}
\begin{theorem}\label{pattern}
{\rm Let ${\bf{b}}$ be given in (\ref{mul}).
Then ${\bf{b}}$ is generalized $k$-regular.}
\end{theorem}
\begin{proof}
{\rm
By the definition of ${\bf{b}}$, for any non-negative integers $e$, $l$ and $j$ with $0\le l\le k^{d-1}-1$ and $0\le j\le k^{e}-1$, we have
\begin{align}\label{produ}
&b(k^{e}(k^{d-1}n+l)+j)=\prod_{y=0}^\infty d_3(k^{e}(k^{d-1}n+l)+j;p_d k^y+\cdots +p_1k^{y+d-1})=\nonumber \\
&\prod_{y=e}^\infty d_3(k^{e}(k^{d-1}n+l)+j;p_d k^y+\cdots +p_1k^{y+d-1})\times \prod_{y=0}^{e-1} d_3(k^{e}(k^{d-1}n+l)+j;p_d k^y+\cdots +p_1k^{y+d-1}).\end{align}
By the definition of $d_3$ and \eqref{reprep}, we have 
\begin{align}\label{produ3}
\prod_{y=e}^\infty d_3(k^{e}(k^{d-1}n+l)+j;p_d k^y+\cdots +p_1k^{y+d-1})=\prod_{y=e}^\infty d_3(k^{e}(k^{d-1}n+l);p_d k^y+\cdots +p_1k^{y+d-1}).\end{align}
For any non-negative integers $e$ and $l$ with $0\le l\le k^{d-1}-1$, we define the sequence ${{\bf{b}}_{e,l}}={(b_{e,l}(n))}_{n\ge0 }$ as $$b_{e,l}(n)=\prod_{y=e}^\infty d_3(k^{e}(k^{d-1}n+l);p_d k^y+\cdots +p_1k^{y+d-1}).$$
By the definition of $d_2$ and \eqref{reprep}, we have
\begin{align}\label{produ4}
\prod_{y=0}^{e-1} d_3(k^{e}(k^{d-1}n+l)+j;p_d k^y+\cdots +p_1k^{y+d-1})=\prod_{y=0}^{e-1} d_3(k^{e}l+j;p_d k^y+\cdots +p_1k^{y+d-1}).\end{align}
For any non-negative integers $e$, $l$ and $j$ with $0\le l\le k^{d-1}-1$ and $0\le j\le k^{e}-1$, we define the constant $d_{e,l,j}$ as $$d_{e,l,j}=\prod_{y=0}^{e-1} d_3(k^{e}l+j;p_d k^y+\cdots +p_1k^{y+d-1}).$$
By \eqref{produ}, \eqref{produ3}, \eqref{produ4}, we get
\begin{align}\label{produ5}
b(k^{e}(k^{d-1}n+l)+j)=d_{e,l,j}b_{e,l}(n).
\end{align}
By \eqref{produ5}, the set $S_{e+d-1}({\bf{b}})$ is contained in the $\mathbb{C}$-subspace generated by the $k^{d-1}$ sequences ${{\bf{b}}_{e,l}}$ ($0\le l\le k^{d-1}-1$). Therefore, the sequence ${\bf{b}}$ is generalized $k$-regular.

\hfill
$\square$}
\end{proof}
From the case $d=1$ in Theorems \ref{pattern1} and \ref{pattern}, the set of generalized $k$-regular sequences includes the set of $k$-additive sequences and $k$-multiplicative sequences. Therefore, the notion of generalized $k$-regular sequences unifies that of $k$-regular, $k$-additive and $k$-multiplicative sequences.
\subsection{Generalized $k$-regular sequences related to Mahler functions}
\quad
In this subsection, we give examples of two classes of generalized $k$-regular sequences related to Mahler functions.
\subsubsection{Generalized $k$-regular sequences related to certain infinite sum of polynomials}
Let $L$ be an integer with $L\ge2$.
Duverney-Nishioka \cite{DuN} investigated the transcendence of values of the following infinite series
\begin{align}\label{DuN}
 \sum_{e=0}^\infty \frac{\sum_{s=1}^{L} a_{s,e} z^{sk^e}}{1+\sum_{s=1}^{L} b_{s,e} z^{sk^e}}
\end{align}
under certain  arithmetical conditions $a_{s,e}, b_{s,e}\in \mathbb{C}$.
Now we consider \eqref{DuN} under the restriction $b_{s,e}=0$.
Let ${\bf{a}}:={(a(n ))}_{n\ge0}$ be a sequence and  $f(z)$ be its generating series defined by 
\begin{align}
f(z)=\sum_{n=0}^{\infty}a(n)z^n:=\sum_{e=0}^\infty (\sum_{s=1}^{L} a_{s,e} z^{sk^e}),\;\;\:a_{s,e}\in \mathbb{C}.
\end{align}
We deduce from Theorem \ref{th2} that the sequence ${\bf{a}}$ is a generalized $k$-regular sequence.
For the proof, we define the following inifinite sum of polynomials
$$f_n(z):=\sum_{e=n}^\infty (\sum_{s=1}^{L} a_{s,e} z^{sk^{e-n}}) \;\;\; (n\ge0 ).$$
As an element of $\mathbb{C}[[z]]$, $f(z)$ satisfies a chain of Mahler type functional equations,
\begin{align}\label{sum}
 f_n(z)=(\sum_{s=1}^{L} a_{s,n} z^{s})+f_{n+1}(z^k) \;\;\; (n\ge0 ).
\end{align}
For any positive integer $l$, let us denote $\mathbb{C}[z]_l$ the set of polynomials in $\mathbb{C}[z]$ whose degrees are less than $l$.
Then, under the expression $L=kS+T$ where $S\ge0$ and  $0\le T\le k-1$, we have
$$\mathbb{C}[z]_{L+1}\subset  \mathbb{C}[z]_{k}+\mathbb{C}[z]_{k}z^k+\mathbb{C}[z]_{k}z^{2k}+\cdots+\mathbb{C}[z]_{k}z^{Sk}.$$
It follows from \eqref{sum} with this observation that
\begin{align*}
&\left(
    \begin{array}{c}
      f_n(z) \\
      1 \\
      z\\
      \vdots \\
      z^S
    \end{array}
  \right)=A_n(z)\left(
    \begin{array}{c}
      f_{n+1}(z^k) \\
      1 \\
      z^k\\
      \vdots \\
      z^{kS}
    \end{array}
  \right).
  \end{align*}
\normalsize
where $A_n(z)$ is an $(S+2)\times (S+2)$ matrix with entries in $\mathbb{C}[z]_{k}$. Hence,
by Theorem \ref{th2}, ${\bf{a}}$ is a generalized $k$-regular sequence. \newline
\newline
\quad
We give the concrete example as follows.
\begin{Exa}
{\rm Recently, Ooto \cite{oo} proved the existence Mahler's $T$-number in positive charecteristic. (See \cite{oo} for definition)
As a particular case he showed that, for a sequence ${(a(n ))}_{n\ge0}$ in ${\mathbb{F}}_{2}$ defined by
\begin{displaymath}
a(n):=\begin{cases} 1 & \mbox{ if $n=2^{4^k l}$ for some integer $k\ge0$}\\
                              &         \mbox{and odd integer \;$l>0$,}\\
                          0 & \mbox{otherwise,}
                                            \end{cases}
\end{displaymath}
the Laurent series $\sum_{n=0}^\infty a(n) T^{-n}$ is a Mahler's $T$-number in ${\mathbb{F}}_{2}((T^{-1}))$. Note that the sequence ${(a(n ))}_{n\ge0}$ is also a generalized $k$-regular sequence in ${\mathbb{F}}_{2}$.}
\end{Exa}
\subsubsection{Generalized $k$-regular sequences related to certain infinite product of polynomials}
Let $L$ be an integer with $L\ge2$. 
Tachiya \cite{Ta}, Amou-V\"a\"an\"anen \cite{AmV} investigated the transcendence of values of the following infinite product
\begin{align}\label{Ta,AmV}
\prod_{e=0}^\infty \frac{(1+\sum_{s=1}^{L} a_{s,e} z^{sk^e})}{ (1+\sum_{s=1}^{L} b_{s,e} z^{sk^e})}
\end{align}
under certain arithmetical conditions on $a_{s,e}, b_{s,e}\in \mathbb{C}$.
As in previous case, we consider \eqref{Ta,AmV} under the restriction $b_{s,e}=0$.
Let ${\bf{a}}={(a(n ))}_{n\ge0}$ be a sequence and $f(z)$ be its generating series defined by
\begin{align*}f(z):=\sum_{n=0}^{\infty}a(n)z^n=\prod_{e=0}^\infty (1+\sum_{s=1}^{L} a_{s,e} z^{sk^e}),\;\;\;a_{s,e}\in \mathbb{C}.
\end{align*}
We deduce from Theorem \ref{th2} that the sequence ${\bf{a}}$ is a generalized $k$-regular sequence.
For the proof, we define the following infinite product of polynomials
$$f_n(z):=\prod_{e=n}^\infty (1+\sum_{s=1}^{L} a_{s,e} z^{sk^{e-n}}) \;\;\; (n\ge0 ).$$
As an element of $\mathbb{C}[[z]]$, $f(z)$ satisfies a chain of Mahler type functional equations
\begin{align}\label{product}
 f_n(z)=(1+\sum_{s=1}^{L} a_{s,n} z^{s})f_{n+1}(z^k).
\end{align}
We choose a positive integer $S$ such that $Sk\le S+L \le Sk+k-1$. Then we can deduce from \eqref{product} together with relation
$$\mathbb{C}[z]_{S+L+1}\subset  \mathbb{C}[z]_{k}+\mathbb{C}[z]_{k}z^k+\mathbb{C}[z]_{k}z^{2k}+\cdots+\mathbb{C}[z]_{k}z^{Sk}$$
that
\begin{align}\label{k-pro}
&\left(
    \begin{array}{c}
      f_n(z) \\
      zf_n(z) \\
      \vdots \\
      z^Sf_n(z)
    \end{array}
  \right)=A_n(z)\left(
    \begin{array}{c}
      f_{n+1}(z^k) \\
      z^kf_{n+1}(z^k) \\
      \vdots \\
      z^{kS}f_{n+1}(z^k)
    \end{array}
  \right),
  \end{align}
where $A_n(z)$ is an $(S+2)\times (S+2)$ matrix in $\mathbb{C}[z]_k$. Hence,
by Theorem \ref{th2}, ${\bf{a}}$ is a generalized $k$-regular sequence.

\section{Acknowledgements}
We would like to thank Professor Masaaki Amou for valuable comments and warm encouragement.
We would like to thank Dr. Tomohiro Ooto for valuable comments. This research was supported by Waseda University Grant for Special Research Project 2017K-178.

\end{document}